\documentclass[11pt,a4paper]{amsart}%\amsart
\usepackage[english]{babel}
\usepackage{amssymb,amsmath,amsthm}
\usepackage{amsfonts}
\usepackage{graphicx}
\usepackage{float}
\usepackage{verbatim}
\usepackage{hyperref}
\hypersetup{colorlinks=true,linkcolor=black,filecolor=magenta,urlcolor=cyan,pdftitle={Overleaf Example},pdfpagemode=FullScreen}
\newtheorem{theorem}{Theorem}
\newtheorem{lemma}[theorem]{Lemma}
\newtheorem{corollary}[theorem]{Corollary}
\newtheorem{proposition}[theorem]{Proposition} 
\theoremstyle{definition} 
 
\newtheorem{remark}[theorem]{Remark}
\newtheorem{hyp}[theorem]{Hypothesis}
\title{Global centers in piecewise linear equations in the cylinder}

\author{J.L. Bravo}
\address{Departamento de Matemáticas, Universidad de Extremadura, 06071 Badajoz, Spain}
\email{trinidad@unex.es}
\author{R. Trinidad-Forte}
\address{Departamento de Ciencias Matemáticas e Informática, Universitat de les
Illes Balears, 07122 Palma, Spain}
\email{
roberto-sebastian.trinidad@uib.cat}

\thanks{
The authors are partially
supported by the project PID2023- 151974NB-I00 funded by MICIU/
AEI/10.13039/ 501100011033/FEDER, UE.}

\begin{document}

\begin{abstract}
We characterize global centers (all solutions are periodic) of the piecewise linear equation $x'=a(t)|x| + b(t)$ when the coefficients $a,b$ are trigonometric polynomials, under some generic hypotheses. 

We prove that the global centers are those determined by the composition condition on $a,b$. That is, the equation has a global center if and only if there exist polynomials $P, Q$ and a trigonometric polynomial $h$ such that $a(t)=P(h(t))h'(t)$,
$b(t)=Q(h(t))h'(t)$.
\end{abstract}
%%%%%%%%%%%%%%%%%%%%%%%%%%%%%%%%%5

\maketitle

\section{Motivations and main result}

Piecewise linear systems are currently being studied due to the accurate description that these systems provide of some phenomena in nature. Indeed, piecewise linear systems seem to present almost the same dynamical behavior as general nonlinear systems (for more information, see \cite{Buzzi,Carmona, Carmona3} and the references therein).
In addition, many non-linear systems can be adequately modeled by continuous piecewise linear systems separated by one or two parallel hyperplanes, splitting the phase plane into different linear dynamic regions (see \cite{Carmona2}).

\medskip

Several authors have studied the scalar case. For instance, in \cite{Carmona}, from a three-dimensional piecewise linear system with two zones, the following reduced one-dimensional $2\pi-$periodic ODE is obtained 
\begin{equation*}
x'=
\begin{cases}
a^+x + b(t),\quad \text{if }x\geq 0,\\
a^-x + b(t),\quad \text{if }x < 0,\\
\end{cases}
\end{equation*}
where $a^+,a^- \in \mathbb{R},\, b(t) = b_0 + b_1\cos(t), \, b_0 \in \{ 0,1 \},$ and $b_1 \geq 0$. The analysis of the limit cycles of this equation determines the dynamics of the three-dimensional piecewise linear system.

\medskip

In this paper, we consider the scalar linear piecewise differential equation
\begin{equation}\label{eq:main}
x' = a(t)|x| + b(t),
\end{equation}  
where $a(t),\,b(t)$ are analytic $\mathcal{C}^1$-functions, $T$-periodic with respect to $t$ and such that $b(t)$ has simple zeroes. Under these hypotheses, a unique maximal solution to the initial value problem exists. A solution $x(t)$ of equation \eqref{eq:main} is periodic if $x(0) = x(T)$. Furthermore, a periodic solution is denoted as \textit{limit cycle} if it is isolated in the set of periodic solutions. We say that a periodic solution is a \textit{center} if it is not a limit cycle, and that a center is \textit{global} if every solution is periodic.

\medskip

In \cite{Gasull}, an upper bound for the number of limit cycles of equation \eqref{eq:main} is obtained by analyzing the displacement map, where $a(t)$ or $b(t)$ are real functions such that one of them has a constant sign.

\medskip 

In contrast, when $a(t)$ or $b(t)$ changes sign in \eqref{eq:main}, no uniform upper bound exists for the number of limit cycles. In fact, in \cite{Coll} (see also \cite{Tineo}) it is proved that for any natural number $k \geq 2$ and $\epsilon$ small enough, the following ODE $x' = \epsilon \cos(kt)|x| + \sin(t)$ has at least $k - 2$ limit cycles. 

\medskip

The center problem studies the necessary and sufficient conditions for an ODE to have a center. The question has been intensively studied for the Abel equations (see e.g. \cite{Alwash, Gine2019}, and references therein). In this paper, the focus is set on the global centers. The statement of the problem goes as follows:
\medskip

\textit{Is it possible to give explicit conditions on $a(t)$ and $b(t)$ such that all the solutions of \eqref{eq:main} are periodic?}
\medskip

Although the number of limit cycles has been studied in several articles (see \cite{Paper2,Huang,Tian} and references therein), the centers of scalar piecewise linear equations have not been profusely studied. Note that there are no equilibrium points and that the return map is not analytic in the entire domain (see~\cite{Tineo}). In particular, there are ODEs with centers formed by solutions with a definite sign that coexist with limit cycles (defined by solutions with no definite sign, see again~\cite{Tineo}).

\medskip

In \cite{Paper2}, the authors gave a characterization for a global center of~\eqref{eq:main} when $a(t),b(t)$ are trigonometric polynomials of degree one. Inspired by this result, we study the general case where the coefficients can be any of any degree. Denote $u(t,x)$ to the solution of~\eqref{eq:main} determined by $u(0,x)=x$.
We assume that the following hypothesis holds, which we prove to be generic in Remark~\ref{re:gen}.

\begin{hyp}\label{hyp}
There exists an open interval $(x_1,x_2)$ such that $u(t,x)$ has exactly two simple zeroes $t_1(x),t_2(x)\in (0,2\pi)$ for any $x\in(x_1,x_2)$, and such that
\[
\lim_{x\to x_1^+} t_2(x)-t_1(x) = 0,\quad \text{or} \quad \lim_{x\to x_2^-} t_2(x)-t_1(x) = 0.
\]
\end{hyp} 

Our main result is the following. 
\begin{theorem}\label{theo:main}
Assume Hypothesis~\ref{hyp} holds. Equation \eqref{eq:main} has a global center if and only if there exist polynomials $P, Q \in \mathbb{R}[x]$ and a trigonometric polynomial $h\in\mathbb{R}[\sin(t),\cos(t)]$, such that 
\begin{equation}\label{eq:cc}
a(t) = P(h(t)) h'(t), \quad 
b(t) = Q(h(t)) h'(t).    
\end{equation} 
\end{theorem}

Equations with centers such that~\eqref{eq:cc} holds are called composition centers. These centers were first introduced by Alwash and Lloyd in \cite{Alwash} and later studied by Briskin, Francoise and Yomdin in \cite{Yomdin, Yomdin2}. Indeed, they conjectured that all centers of the Abel equation are composition centers. Although the conjecture was proven to be false for $a,b$ trigonometric polynomials or even polynomials (\cite{Pakovich2002, Gine2019}), they seem to play an important role in the set of centers for the Abel equations (\cite{Yomdin2,Pakovich3}). By Theorem~\ref{theo:main}, they are also relevant in the case of piecewise-linear equations. 

\medskip

The paper is organized as follows. In Section 2, we include some preliminary results and develop some elementary ideas; Section 3 contains a subsection about the solutions and their structure, and the final section contains the proof of the main result.

\section{Preliminary results.}

This section contains some preliminary results and tools that will help us study periodic solutions and establish an adequate basis to deal with the functions and polynomials involved.

\subsection{Zero-dimensional Abelian integrals}

Assume that Hypothesis~\ref{hyp} holds. From this condition, we will obtain information about the coefficients $a(t)$ and $b(t)$ of \eqref{eq:main}. To do so, we need some results provided as a consequence of Lüroth's theorem. 

\medskip 

%These were developed by Álvarez, Bravo and Marde\v si\'c in \cite{Amelia} and are included in the paper for reference.

Let $f \in \mathbb{C}(z)$ be a rational function of degree $m > 1$. Let $(z_i(c)))_{1\leq i\leq m}$ denote an $m-$tuple of analytical preimages $z_i(c) \in f^{-1}(c)$, where $c \in \mathbb{C} \setminus \Sigma$ and $\Sigma$ is the set of critical values of $f$. Define a \textit{zero-dimensional cycle} (shorter cycle) $C(t)$ of $f$ as the sum
\begin{equation*}
    C(t) = \sum_{i=1}^mn_iz_i(t) \quad \text{with} \quad \sum_{i=1}^m n_i = 0, \quad  n_i \in \mathbb{Z}. 
\end{equation*}

We say that a cycle $C(t)$ is \textit{simple} if it is of the form $C(t) = z_j(t) - z_i(t)$, and \textit{trivial} if $n_i = 0$ for every $i$.

\medskip 

Let $g \in \mathbb{C}(z)$ be a rational function. As in \cite{Amelia}, define \textit{zero-dimensional Abelian integrals of $g$ along the cycle $C(t)$} as
\begin{equation*}
    \int_{C(t)}g := \sum_{i=1}^m n_i g\left( z_i(t)\right).
\end{equation*}

\begin{proposition}(see Proposition 5.3 from \cite{Amelia})\,\label{prop:amelia}
    Let $f(z) \in \mathbb{C}(z)$ and consider a simple cycle $C(t) = z_i(t) - z_j(t).$ Then
    \[
    \int_{C(t)} g \equiv 0
    \]
    if and only if there exist $f_0,g_0,h \in \mathbb{C}(z)$ such that $f = f_0 \circ h, g = g_0 \circ h$ and $h(z_i(t)) = h(z_j(t))$.
\end{proposition}

\subsection{Laurent and Trigonometric Polynomials.}
Next, we cover some results about the relationship between trigonometric and Laurent polynomials. This is done since some of the applied results require a Laurent polynomial, which we will later use for our trigonometric context. This consideration of the problem changes the domain from $[0,2\pi]$ to the complex unit circle. 
\medskip

Let $p(t) \in \mathbb{R}[\sin\,t,\cos\,t]$ be a real trigonometric polynomial of degree $n$. Then, it can be written as
\[
p(t) = a_0 + \sum_{k=1}^n a_k \cos(k t) + b_k \sin(k t), 
\]
where $a_0,a_k,b_k \in \mathbb{R},\, \forall k \in \{1,...,n \}.$
Let us denote by $\mathcal{L}=\mathbb{C}[z,z^{-1}]$ the ring of Laurent polynomials with complex coefficients. Define the following sub-ring,
\[
\mathcal{A} = \{p(z) \in \mathbb{C}[z,z^{-1}]: p(z) = \sum_{k=0}^n \omega_k \,  z^k + \dfrac{\bar{\omega}_k}{z^k},\, \omega_k \in \mathbb{C},\, n\in \mathbb{N}\}.
\]
A Laurent polynomial $p \in \mathcal{A}$ has degree $n$ if 
\[
p(z) = \sum_{k=0}^n \omega_k \,  z^k + \dfrac{\bar{\omega}_k}{z^k},\quad \text{with }\omega_n\neq 0.
\]
Note that the definition does not depend on the imaginary part of $\omega_0$.

\begin{proposition}\label{prop:laurent}
There is a bijective correspondence between trigonometric polynomials with real coefficients of degree $n$, and the set $\mathcal{A}$ of Laurent polynomials of degree $n$.  

\medskip 

Concretely, if $p\in\mathcal{A}$, then 
$p(e^{it})\in\mathbb{R}[\sin\,t,\cos\,t]$.
\end{proposition}
\begin{proof}
Let $p(t)$ be a trigonometric polynomial. Recall that it can be written as 
\[p(t) = a_0+ \sum_{k=1}^n a_k \cos(k\, t) + b_k \sin(k\, t).\] Through the change of variables $z=e^{it}$, \[\cos(t) = \frac{z + z^{-1}}{2},\quad \sin(t) = \frac{z - z^{-1}}{2i},\]
so we obtain the Laurent polynomial $p\in\mathcal{A}$ defined by
\[p(z) = a_0 + \frac{1}{2} \sum_{k=1}^n (a_k -  b_k\, i)z^k + (a_k + b_k\, i) \dfrac{1}{z^k}.  \]
    
The transformation from a trigonometric polynomial to a Laurent one is injective. The proof follows by taking into account the dimensions as real vector spaces.
\end{proof}

The next lemma gives us more information on the factorization of the functions. We include an adaptation of the proof given by Zieve in \cite[Lemma 2.1]{Zieve} as we will need some details on the construction of the factorization. For more details, we refer the reader to \cite[p. 25]{Pakovich} and \cite{Pakovich2}.

\begin{lemma}\label{lemma:1} 
If $p = g \circ h$, where $p \in \mathcal{A}$ and $g,h \in \mathbb{C}(z)$, then there is a degree-one $\mu \in \mathbb{C}(z)$ such that $G := g \circ \mu$ and $H := \mu^{-1} \circ h$ satisfy one of the following cases:
    \begin{itemize}
        \item $G \in \mathbb{R}[x]$ and $ H \in \mathcal{A}$;
        \item $G \in \mathcal{A}$ and $H = z^k$ for some $k \in \mathbb{N}$.
    \end{itemize}
\end{lemma}
\begin{proof}
The poles of $p = g\circ h$ are the preimages under $h$ of the poles of $g$. By hypothesis, the poles of $p$ are $\{0,\infty \}$. Hence, $g$ has at most two poles.
    
Let us suppose that $g$ has a unique pole, $\alpha\in\mathbb{C}$. Hence, $h(0) = h(\infty) = \alpha$. Let a degree-one $\mu \in \mathbb{C}(x)$ for which $\mu (\infty) = \alpha$, so that $G := g \circ \mu$ has $\infty$ has its unique pole. This implies $G \in \mathbb{C}[z]$. Then, $p = G \circ H$, where $H := \mu^{-1} \circ h$. Note that $H$ has $\{0,\infty\}$ as poles, so $H \in \mathcal{L}$. 

To conclude, we need to show that it is possible to choose $G \in \mathbb{R}[x]$ and $ H \in \mathcal{A}$. To do so, we construct $\mu$ y $\mu^{-1}$ explicitly.
%First, \color{red} since $g$ has a unique pole, say $\alpha$, then $g(z)$ is a polynomial on $\frac{1}{z-\alpha}$. \color{black}\marginpar{\ro{Tengo mis dudas sobre si esto en rojo es necesario}}
Since we have taken $\mu(\infty) = \alpha$, we have \[\mu(x) = \dfrac{\alpha x + b}{x + d},\quad \mu^{-1}(x) = \dfrac{b - dx}{x - \alpha},\] where $b,d\in \mathbb{R}$. Note that  both maps depend on two parameters, which can be arbitrary chosen. Therefore, we choose them so that $H := \mu^{-1}\circ h \in \mathcal{L}$ has no independent term, and such that its principal coefficients (highest and lower degree) are complex-conjugate.

If we now take the coefficients of $G \circ H$ and compare them in pairs with the coefficients of $p$, starting with the highest and the lowest degrees, in a finite number of steps we obtain that $G$ is a real polynomial and that $H \in\mathcal{A}$.

%For a better understanding, we include the two first examples. Let $G(z) = \sum_{i=0}^n g_iz^i,\, g_i \in \mathbb{C},\, i=0,...,n$, and let $H(z) = \sum_{i=0}^m h_{-m}z^{-m} + h_mz^m,\, h_i \in \mathbb{C}, \, i=-m,...,m$, such that $h_m = \overline {h_{-m}}$ and $h_0 = 0$. Now, equalling the coefficients of $p$ and $G\circ H$ of highest and lowest degree, we get that $g_nh_m^n = \overline {g_nh_{-m}^n}$. Therefore, $g_n \in \mathbb{R}$. Lowering the degree of the coefficients, we get that $ng_nh_{m-1}h_m^{n-1}= \overline {ng_nh_{-m+1}h_{-m}^{n-1}}$, and so we get that $h_{m-1} = \overline{ h_{-m+1}}$. If we continue with the process, it follows that $G \in \mathbb{R}[x]$ and $ H \in \mathcal{A}$.
\medskip

Suppose that $g$ has two poles, say $\alpha$ and $\beta$. Since $g \circ h$ has at most two poles, both $\alpha$ and $\beta$ must have unique preimages under $h$, which must be $0$ and $\infty$. Suppose that $\alpha = h(0)$, $\beta = h(\infty)$ and $\gamma = h(1)$. Let us take a degree-one $\mu \in \mathbb{C}(z)$ such that $\mu(0) = \alpha$, $\mu(\infty) = \beta$ and $\mu(1) = \gamma$. Then the poles of $G := g \circ \mu$ are $0$ and $\infty$, so $G \in \mathcal{L}$, and $H := \mu^{-1} \circ h$ has its unique pole at $\infty$, so $H \in \mathbb{C}[z]$. Note that $0$ is the only root of $H$, so $H$ is a monomial such that $H(1) = 1$. Therefore, $H$ is monic. Now, as $p\in\mathcal{A}$, it is easy to check that $G\in\mathcal{A}$.
\end{proof}

%\begin{remark}
%Since $\mathcal{A} \subset \mathcal{L}\setminus \left(\mathbb{C}[z]\cup\mathbb{C}[z^{-1}]\right)$, the same result is obtained if, in particular, $p \in \mathcal{A}$. 
Due to the correspondence given by Proposition \ref{prop:laurent}, it is normal to ask if there is a similar factorization when $p$ is a trigonometric polynomial.
%\end{remark}
%Proposition~\ref{prop:laurent} establishes a bijective correspondence between trigonometric polynomials with real coefficients and  Laurent polynomials in $\mathcal{A}$. 
To simplify the exposition, if $p(t)$ is a trigonometric polynomial, we denote the corresponding Laurent polynomial as $p(z)$.

\medskip

\begin{corollary}\label{corol:1}
Let $p(t)$ be a trigonometric polynomial with real coefficients. If $p(z) = g(h(z))$ for every $z\in\mathbb{C}$, where $g,h \in \mathbb{C}(z)$, then $p(t) = G \circ H$, where $G,H$ satisfy one of the following cases:
    \begin{itemize}
        \item $G \in \mathbb{R}[x]$ and $ H\in \mathbb{R}[\sin\,t,\cos\,t]$;% is trigonometric with real coefficients;
        \item $G\in \mathbb{R}[\sin\,t,\cos\,t]$ and $H = k t$ for some $k \in \mathbb{N}$.
    \end{itemize}
\end{corollary}
\begin{proof}
Since $p(z) = g \circ h$, the conditions for Lemma \ref{lemma:1} hold, and we get a new  factorization $p(z) = G \circ H$, such that $G := g \circ \mu$ and $H := \mu^{-1} \circ h$, with $\mu \in \mathbb{C}(z)$.
Let us now study the two possible cases from Lemma~\ref{lemma:1} separately.

\medskip

Assume that there exist $G \in \mathbb{R}[x]$, and $H \in \mathcal{A}$ such that $p(z)=G(H(z))$. Then 
$p(t) = G(H(e^{it}))$, so defining $H(t) = H(e^{it})\in \mathbb{R}[\sin\,t,\cos\,t]$, $p(t)=G\circ H$, where $G$ is a real polynomial and $H$ is a real trigonometric polynomial. 
\medskip 

Assume that there exist $G \in \mathcal{A}$ and $H = z^k$ for some $k \in \mathbb{N}$ such that $p(z) = G(H(z))$. Then 
\[
p(t)=G(H(e^{it}))=G(e^{i k t}).
\]
Taking $G(t) = G(e^{it})\in \mathbb{R}[\sin\,t,\cos\,t]$, the result follows. 
\end{proof}

\section{Structure of the set of solutions}

In this section, we study the signs of the solutions and the conditions for periodicity. The section is organized as follows: First, we analyze the structure of the solutions in the general setting. Second, we obtain the conditions for global centers and some preliminary results for the coefficients. Finally, we state and prove some useful results to carry out the proof of the main result.

\medskip 

Consider the general piecewise equation 
\begin{equation}\label{eq:gen}
x'= h(t,x)=
\begin{cases}
f(t,x),\quad \text{if }x\geq 0,\\
g(t,x),\quad \text{if }x < 0,\\
\end{cases}
\end{equation}
where $f(t,x)$ and $g(t,x)$ are continuous and locally Lipschitz continuous with respect to $x$ functions, $T$-periodic with respect to $t$,  and $f(t,0)=g(t,0)$. The structure of the zeroes of the solutions of~\eqref{eq:gen} is described in \cite{Tineo}, where the following result is proved and which allows us to separate the plane into bands of solutions, each of them with the same number of zeroes within the band.

\begin{proposition}[\cite{Tineo}]\label{prop:3}
Let $f,g:\mathbb{R}^2\to\mathbb{R}$ be continuous, locally
Lipschitz continuous with respect to $x$, and $T$-periodic with
respect to $t$ functions. Suppose that $c(t) :=f(t,0) = g(t,0)$ has exactly $n$
zeroes in $[0,T]$ and that the solutions of (\ref{eq:gen}) are
defined for every $t\in\mathbb{R}$. Then there exist
$$-\infty=x_0<x_1<\ldots<x_r<x_{r+1}=+\infty,\quad 1\leq r\leq
n+2,$$ such that, for each $i\in\{0,\ldots,r\}$ and
$x\in(x_i,x_{i+1})$, the zeroes of the map
$$t\in[0,T]\to u(t,0,x),$$ are simple. Moreover, for each $x\in(x_i,x_{i+1})$, their number
is constant and bounded by $n+1$.
\end{proposition}

Consider now \eqref{eq:main} with $a,b$ trigonometric polynomials. It satisfies the conditions of Proposition~\ref{prop:3} as the equation is globally Lipschitz with respect to $x$, so every solution is defined for every $t\in\mathbb{R}$, and $c(t) = b(t)$ has at most $2m$ roots, where $m$ is the degree of $b$. 

Moreover, by a translation in time, we assume that $b(0) = 0$. In particular, there exists $i\in\{1,\ldots,r\}$ such that $x_i=0$. Arguing as in~\cite{Tineo}, one obtains that $x_1<x_2<\ldots<x_r$ can be chosen as the value at $t=0$ of the solutions with a non-simple zero.

\medskip

\begin{corollary}\label{coro:definite_sign}
Let $u(t)$ be a periodic solution such that $u(0) \neq x_1,\ldots,x_r$. Then $u$ has an even number
of zeroes, all simple. 

Moreover:
\begin{enumerate}
    \item If $u(0)>x_r$, then $u(t)>0$ for every $t\in\mathbb{R}$ and if $u(0)<x_1$, then $u(t)<0$ for every $t\in\mathbb{R}$.
    \item Assume that $x_1,\ldots,x_r$ are taken as the value at $t=0$ of the solutions with a non-simple zero. If $x_r<u(0)<x_1$, then $u$ has non-definite sign (it is strictly positive for some values of $t$ and strictly negative for some other values).
\end{enumerate}

\end{corollary}
\begin{proof}
Recall that $x_1,\ldots,x_r$ are chosen as the value at $t=0$ of the solutions with a non-simple zero. Then, by the uniqueness of solutions, a periodic solution $u(t)$ with $u(0) \neq x_1,...,x_r$ must have only simple zeroes. There are an even number of them, otherwise it would not be periodic.

According to Proposition \ref{prop:3}, the solution $u(t,x_r)$ such that $u(0,x_r)=x_r$ is the border between bands of solutions with different number of zeroes. In particular, $u(t,x_r)$ has one non-simple zero by hypothesis and it verifies $u(t,x_r)\geq 0 \, ~\forall t \in [0,2\pi)$. As a consequence, if $u(0) > x_r$, then $u(t) > 0$ for all $t \in \mathbb{R}$ due to the uniqueness of the solutions. The negative case is analogous.

If $u(t)$ has definite sign, for instance positive, and $x_1 < u(0) < x_r$, then all the solutions above $u(t)$ must have definite sign and no zeroes, including the one that takes the value $x_r$ at $t=0$. This is a contradiction, since $u(t,x_r)$ has a non-simple zero. The same conclusion is reached if $u(t)$ is negative.

\end{proof}

Let $u(t,\bar x)$ denote the solution of~\eqref{eq:main} determined by the initial condition $u(0,\bar x)=\bar x$. Let $t_1<t_2<\ldots<t_k\in(0,2\pi)$ be zeroes of $t\to u(t,\bar x)$. 

\begin{remark}\label{re:gen}
If $\bar x=x_r$, by Corollary~\ref{coro:definite_sign} and continuity, $u(t,x_r)\geq 0$. Note that if $u(t,x_r)$ is strictly positive, then $x_r$ can be removed from the values $\{x_1,\ldots,x_r\}$ and Proposition~\ref{prop:3} still holds. Therefore, there is no restriction assuming that $u(t,x_r)$ has a zero, which is necessarily non-simple. Moreover, $t_1,t_2\in(0,2\pi)$ are two zeroes of $u(t,x_r)$ if and only if $t_1,t_2$ are zeroes of $b$ and 
\[
\int_{t_1}^{t_2} b(t)\exp \left( \int_t^{t_2}a(s)ds \right) dt = 0.
\]
That is, generically, $u(t,x_r)$ has exactly one zero. Moreover, that implies that $u(t,x)$ has exactly two zeroes in $(0,2\pi)$ for $x<x_r$ close enough and that Hypothesis~\ref{hyp} holds. 
\end{remark}

%\begin{figure}[h]
%    \centering
%    \includegraphics[scale=0.3]{esquema.pdf}
%    \caption{\ro{Tengo que añadir en la imagen los $x_i$ con Latex, me da problemas el Inkscape.}}
%\end{figure}

\subsection{Solutions with two simple zeroes}

In the following, we assume that there exists $\bar x\not\in\{ x_1,\ldots,x_r \}$ such that $u(t,\bar x)$ has exactly two simple zeroes in $(0,2\pi)$. By Proposition~\ref{prop:3}, the solution $u(t,x)$ has exactly two simple zeroes in $(0,2\pi)$ for $x$ close to $\bar x$. Recall that $\bar x\neq 0$, as $t=0$ is a zero of $b(t)$, so $x_i=0$ for some $i$. For simplicity, we will assume $\bar x<0$, being the case $\bar x>0$ analogous. 
\medskip 

Define the functions
$t_i(x)$, $i=1,2$ as the $i$-th zero of $t\to u(t,x)$, $t\in(0,2\pi)$, for $x$ in an open interval containing $\bar x$, $I$. By continuous dependence of the solutions with respect to the initial conditions, $t_1,t_2$ depend continuously on the initial conditions. Moreover,

\begin{proposition}\label{prop:4}
Assume $x<0$. Then
\begin{enumerate}
\item $t_1,t_2$ are analytic with analytic inverse. 
\item $b(t_1(x))>0$, $t_1'(x) < 0$, $b(t_2(x))<0$, $t_2'(x)>0$ for every $x\in I$.
\item $t_1,t_2$ are related by the next implicit equation 
\begin{equation}\label{eq:t1t2}
  \int_{t_1(x)}^{t_2(x)} b(t)\exp \left( \int_t^{t_2(x)}a(s)ds \right) dt = 0.    
\end{equation} 
\item For any $x\in I$, $u(t,x)$ is periodic if and only if
\begin{equation}\label{eq:cperiodica}
\int_{t_2(x)}^{t_1(x) + 2\pi} b(t)\exp \left( \int_t^{t_1(x) + 2\pi}-a(s)ds \right) dt = 0.
\end{equation}
\end{enumerate}
\end{proposition}
\begin{proof}

\textit{(1):} Let us fix $\bar{x} \in \mathbb{R}$, and let $\bar{t}_1 \in (0,2\pi)$ be the first zero of the solution. Since the zeroes are supposed simple, $u_t(\bar{t}_1,\bar{x}) \neq 0$, where the subscript denoted partial derivative. Therefore, the conditions for the Implicit Function Theorem hold, and we can apply it to the map $(t,x) \rightarrow u(t,x)$ at $(\bar{t}_1,\bar{x})$. As a consequence, one obtains the existence of an analytic function $t_1(x)$ defined in a neighborhood of $\bar{x}$, such that $u(t_1(\bar{x}),\bar{x}) = 0$ and $t_1(\bar{x}) = \bar{t}_1$.
The process is analogous for the function $t_2(x)$. The existence of the inverse is proved in the following point.

\medskip 

\textit{(2):} As we assumed $\bar x<0$, the solution changes from negative to positive at $t_1$. As it is a simple zero, then $u_t(t_1(x),x) = b(t_1(x)) > 0$. Analogously, since the solution changes from positive to negative at $t_2$, then $u'(t_2(x),x) = b(t_2(x)) < 0$.
    \medskip

Let us now take the derivative of $u(t_1(x),x) = 0$,
\[
u_t(t_1(x),x)t'_1(x) + u_x(t_1(x),x) = 0,
\]
and so
\[
t_1'(x) = - \dfrac{u_x(t_1(x),x)}{u_t(t_1(x),x)}<0,
\]
as $u_t(t_1(x),x)>0$, and
%\[u_x(t,x)= \exp \left( \int_0^{t}-a(s)\,ds\right)>0.\]
%We can deduce the sign of $t_1(x)$ from the expression of the solution, $u(t,x) = x \exp \left( \int_0^{t}-a(s)\,ds\right) + \int_0^{t}b(t) \exp \left(\int_0^{t}-a(s)\,ds \right)\,dt$. 
%If one takes the derivative with respect to $x$ at $t=t_1(x)$, one gets 
    \[
            u_x(t_1(x),x) = \exp \left( \int_0^{t_1(x)}-a(s)\,ds\right) > 0.
    \]
%Taking now the derivative with respect to $t$, we get that $u_t(t_1(x),x) > 0,$ since the solution changes from negative to positive at $t = t_1(x)$.
%Therefore, $t_1'(x) < 0$. 
The process is analogous for $t_2(x)$.

\medskip

Since $a(t)$ and $b(t)$ in \eqref{eq:main} are analytic and also the solution $u(t,x)$ is analytic, so are $t_1(x)$ and $t_2(x)$. As $t_1'(x),\,t_2'(x) \neq 0$, then, we can use the Inverse Function Theorem to obtain an analytic inverse for both $t_1(x)$ and $t_2(x)$.
\medskip

\textit{(3):} Denote $u(t,s,x)$ as the solution such that $u(s,s,x)=x$. Since $t_1(x)$ and $t_2(x)$ are the two zeroes of the solution, it holds \[u(t,t_1(x),0)=u(t,t_2(x),0) = u(t,x).\] Then $u(t,t_1(x),0)$ is non-negative for $t \in [t_1(x),t_2(x)]$. From the expression of the solution of a linear ODE we get
\[
u(t_2(x),t_1(x),0) = \int_{t_1(x)}^{t_2(x)} b(t)\exp \left( \int_t^{t_2(x)}a(s)ds \right) dt = 0.
\]

\textit{(4):} Suppose that $u(t,x)$ is periodic. %Since the field of \eqref{eq:main} is $2\pi$-periodic, it follows that the solution is also $2\pi$-periodic, \textit{i.e.} $u(t,0,x) = u(t+2\pi,0,x)$. 
If $t_1(x)$ is the first zero of the solution, then $t_1(x)+2\pi$ is the first zero of the solution after one period. Therefore, the solution is negative for $t \in (t_2(x),t_1(x)+2\pi)$. From the expression of the solution, we get
\[
u(t_1(x)+2\pi,t_2(x),0) = \int_{t_2(x)}^{t_1(x) + 2\pi} b(t)\exp \left( \int_t^{t_1(x) + 2\pi}-a(s)ds \right) dt = 0
\]
For the converse, suppose that \eqref{eq:cperiodica} holds. Then $t_1(x)+2\pi$ is a zero of the solution. Since $t_1(x)$ is also a zero, the solution is $2\pi$-periodic.
\end{proof}

If in addition the solution $u(t,x)$ is periodic for every initial condition in the interval $I$ where $u(t,x)$ has two simple zeroes, then $a$ has null integral between the two zeroes. 

\begin{lemma}\label{le:2}
If $u(t,x)$ is periodic for every $x\in I$, then
\begin{equation*}%\label{eq:cs}
\int_{t_1(x)}^{t_2(x)} a(s)\,ds =0.
\end{equation*}
\end{lemma}
\begin{proof}
For clarity, we omit the arguments of the functions $t_1$ and $t_2$.    

Assume that $u(t,x)$ is periodic for every $x\in I$. Deriving in \eqref{eq:t1t2} and \eqref{eq:cperiodica} with respect to $x$, we get
\begin{equation}\label{eq:der}
\begin{split}
b(t_2)t_2' - b(t_1)t_1' \exp \left( \int_{t_1}^{t_2}a(s)ds \right) =& \,0,\\
b(t_1)t_1' - b(t_2)t_2' \exp \left( \int_{t_1}^{t_2}a(s)ds \right) =& \, 0,
\end{split}
\end{equation}    

By Proposition~\ref{prop:4}, $b(t_1)t_1',b(t_2)t_2'< 0$. As a consequence,  
\[\int_{t_1}^{t_2} a(s)\,ds =0.\]
\end{proof}

\section{Proof of the main result}

This section is devoted to the proof of Theorem~\ref{theo:main}. 

\begin{proof}
We assume that Hypothesis~\ref{hyp} holds, that is, $u(t,x)$ has two simple zeroes, $t_1(x),t_2(x)\in (0,2\pi)$ for every $x\in I$, where $I$ is an open interval, and such that
\[
\lim_{x\to \min(I)^+} t_2(x)-t_1(x) = 0,\quad \text{or} \quad \lim_{x\to \max(I)^-} t_2(x)-t_1(x) = 0.
\]
For simplicity, assume that the first case holds and let $\bar x  = \min(I)$.  By Proposition~\ref{prop:4}, $t_1,t_2$ are analytic monotonous functions, so the limit of $t_1,t_2$ when $x\to\bar x^+$ exists (and it is a zero of $b$). That is, 
\[\bar t:=\lim_{x\to \bar x^+} t_2(x)=\lim_{x\to \bar x^+} t_1(x),\quad b(\bar t)=0.\]

\medskip

Assume that \eqref{eq:main} has a global center. 
Denote
\[
A(t) =\int_0^t a(s)\,ds,
\quad B(t) =\int_0^t b(s)\,ds.
\]
We will omit the arguments of the functions $t_1,t_2$ for simplicity.

\medskip

By Lemma~\ref{le:2},
\[
A(t_2)-A(t_1) = \int_{t_1}^{t_2} a(t)\,dt = 0.
\]

\medskip

\textit{Claim 1: $B(t_1) = B(t_2)$.}
\medskip

As $\int_{t_1}^{t_2} a(t)\,dt =0$, from \eqref{eq:der}
we obtain
\[
b(t_1)t_1' -b(t_2) t_2' = 0.
\]
Integrating this expression in $(\bar x,x)$, 
\[
B(t_1)-B(t_2) = c,
\]
with $c \in \mathbb{R}$. Note that $\lim_{x \rightarrow \bar{x}} t_1(x) = \bar{t} = \lim_{x \rightarrow \bar{x}} t_2(x)$. Therefore,
\[
\lim_{x \rightarrow \bar{x}}  B(t_1(x)) - B(t_2(x))  = 0,
\]
so $c = 0$, and $B(t_1)=B(t_2)$.

\medskip

\textit{Claim 2: there exist $P,Q \in \mathbb{R}[x]$ and a trigonometric polynomial $H$, such that $A = P \circ H,$ $B = Q \circ H$.}
\medskip

In order to continue, we need to consider $A$ and $B$ as Laurent polynomials. Since $A$ and $B$ are trigonometric, by Proposition \ref{prop:laurent} we can take their Laurent expression $A(z),\,B(z) \in \mathcal{A}$.

As the function $t_1,t_2$ are locally invertible, they can be parametrized as two families of preimages of $A$. Abusing notation, let $t_1,t_2$ be the corresponding preimages of $A(z)$, parametrized by $A(t_1(c))=A(t_2(c))=c$ and let us define the simple cycle $C(c) = t_2(c) - t_1(c)$. In order to use Proposition \ref{prop:amelia}, define the zero-dimensional Abelian integral
\[
\int_{C(c)}B(z) := B(t_2(c)) - B(t_1(c)).
\]
According to Claim 1, $B(t_2(c)) = B(t_1(c))$, and so $\int_{C(c)}B(z) = 0$. By Proposition \ref{prop:amelia}, there exist $A_0,B_0,h \in \mathbb{C}(z)$ such that $A(z) = A_0 \circ h,$ $B(z) = B_0 \circ h$, and $h(t_2(c)) = h(t_1(c))$.
\medskip

After this first factorization of $A(z)$ and $B(z)$, we can use Corollary \ref{corol:1} to obtain a second factorization for their trigonometric equivalent expressions. As a consequence, $A(t) = P \circ H$ and $B(t) = Q \circ H$, where one of the following conditions holds true:
\begin{itemize}
    \item $P,Q \in \mathbb{R}[x]$ and $ H(t)\in\mathbb{R}[\sin\,t,\cos\,t]$;
    \item $P,Q\in\mathbb{R}[\sin\,t,\cos\,t]$ and $H(t) = kt$ for some $k \in \mathbb{N}$.
\end{itemize}
Let us realize that we use Corollary \ref{corol:1} for each $A(z)$ and $B(z)$, which gives us one factorization for each function. However, the $H$ that appears in both of them is unique, since $H$ is determined by the function $h$ from the first factorization $A(z) = A_0 \circ h,$ $B(z) = B_0 \circ h$ (see proof of Lemma \ref{lemma:1}). 

Suppose that the second condition is true. Then, $A,B$ are $2\pi/k$-periodic. Therefore, the same holds for $a,b$ and for the solutions of~\eqref{eq:main}, in contradiction to $u(t,x)$ being a periodic solution with exactly two simple zeroes. 
\medskip

We conclude that there exist $P,Q \in \mathbb{R}[x]$, and a trigonometric polynomial with real coefficients, $H$, such that $A(t) = P \circ H,$ $B(t) = Q \circ H$, or equivalently, $a(t) = P'(H(t))H'(t)$ and $b(t) = Q'(H(t))H'(t)$, with $P',Q'\in \mathbb{R}[x]$ and $H'(t)$ trigonometric.
\medskip

For the inverse implication, a change of variables is required. Denote $\omega = h(t)$, so that $t = h^{-1}(\omega)$, and define $y(\omega) = x(h^{-1}(\omega)) = x(t)$. Therefore,
    \[
    y'(\omega)%=\dfrac{\partial y}{\partial \omega} = \dfrac{\partial x(h^{-1}(\omega))}{\partial \omega} 
    = \dfrac{\partial x(h^{-1}(\omega))}{\partial t} \dfrac{\partial h^{-1}(\omega)}{\partial \omega}
%     \[= \left( P'(h(t))h'(t)|x(t)| + Q'(h(t))h'(t) \right)\dfrac{1}{h'(t)},
%    \]
%    and operating, we get
%    \[
    %y'(\omega) = 
    %\dfrac{\partial y}{\partial \omega} 
    = P'(\omega)|y(\omega)| + Q'(\omega).
    \]
    
    Let us note that $x(0) = y (h(0)) = y(h(2\pi)) = x(2\pi)$, and therefore equation \eqref{eq:main} has a global center.
\end{proof}

\end{document}